\documentclass[graybox]{svmult}
\usepackage{amsmath}
\usepackage{eucal}
\usepackage[active]{srcltx}
\usepackage{amssymb}
 \usepackage[usenames,dvipsnames]{pstricks}
 \usepackage{epsfig}
\usepackage[ansinew]{inputenc} %para acentos em portugues

\newcommand{\p}{\partial}
\newcommand{\pfrac}[2]{\genfrac{}{}{}{1}{#1}{#2}}

\newcommand{\at}[2]{\genfrac{}{}{0pt}{}{#1}{#2}}

\newcommand{\bb}[1]{{\mathbb #1}}
\newcommand{\mc}[1]{{\mathcal #1}}
 \newcommand{\A}{C^{1,2}([0,T]\times [0,1])}

\newcommand{\Y}{\mathcal{Y}}
\usepackage{helvet}         % selects Helvetica as sans-serif font
\usepackage{courier}        % selects Courier as typewriter font

\usepackage{graphicx}        % standard LaTeX graphics tool
\usepackage{multicol}        % used for the two-column index

\begin{document}

\title*{Slowed exclusion process: hydrodynamics, fluctuations and phase transitions}
\author{Tertuliano Franco, Patr\' icia Gon\c calves and Adriana Neumann}
\institute{Tertuliano Franco \at UFBA,
 Instituto de Matem\'atica, Campus de Ondina, Av. Adhemar de Barros, S/N. CEP 40170-110,
Salvador, Brazil.
\\\email{tertu@impa.br}
\and Patrícia Gonçalves \at PUC-RIO, Departamento de Matem\'atica, Rua Marqu\^es de S\~ao Vicente, no. 225, 22453-900, Gávea, Rio de Janeiro, Brazil \\and\\
CMAT, Centro de Matem\'atica da Universidade do Minho, Campus de Gualtar, 4710-057 Braga, Portugal.
\\\email{patricia@mat.puc-rio.br}
\and Adriana Neumann \at UFRGS, Instituto de Matem\'atica, Campus do Vale, Av. Bento Gon\c calves, 9500. CEP 91509-900, Porto Alegre, Brazil.
\\\email{aneumann@impa.br}}

\maketitle

\abstract*{This is a short survey on recent results obtained by the authors on dynamical phase transitions of interacting particle systems. We consider particle systems with exclusion dynamics, but it is conjectured that our results should hold for a general class of particle systems. The parameter giving rise to the phase transition is the ``slowness" of a single bond in the discrete lattice.
The phase transition is verified not only in the hydrodynamics, but also in the fluctuations of the density, the current and the tagged particle. Moreover, we found a phase transition in the \textit{continuum}, that is, at the level of the hydrodynamic equations, in agreement with the dynamical phase transition for the particle systems.}

\abstract{This is a short survey on recent results obtained by the authors on dynamical phase transitions of interacting particle systems. We consider particle systems with exclusion dynamics, but it is conjectured that our results should hold for a general class of particle systems. The parameter giving rise to the phase transition is the ``slowness" of a single bond in the discrete lattice.
The phase transition is verified not only in the hydrodynamics, but also in the fluctuations of the density, the current and the tagged particle. Moreover, we found a phase transition in the \textit{continuum}, that is, at the level of the hydrodynamic equations, in agreement with the dynamical phase transition for the particle systems.}

\section{Introduction}

 A major question in Statistical Mechanics is how  to perform the limit from the discrete to the continuum in such a way that the discretization of the system really gives the correct description of the continuum? This question gave rise to plenty of famous models and results, both in Physics and Mathematics. In the particular context of \emph{particle systems} and \emph{hydrodynamic limits}, the passage of the discrete to the continuum is a consequence of rescaling both time and space. The discrete system consists in a collection of  particles with a stochastic dynamics. Depending on the prescribed  interaction we are lead to different limits.  Therefore the random interaction of the microscopic system is connected to the macroscopic phenomena to be explored.

As the main reference on the subject, we cite the classical book \cite{kl}, which treats the limit of several particle systems, as the zero range process, the symmetric and asymmetric exclusion process, the generalized $K$-exclusion process, independent random walks and some of their scaling limits.
We point out some of the possible natures of those scaling limits.

The  scaling limit for the time-trajectory of the spatial density of particles is the so-called \emph{hydrodynamic limit} of the system, which is a Law of Large Numbers (L.L.N.) type-theorem. The scaling limit for how the discrete system oscillates around its hydrodynamic  limit is usually referred as  \emph{fluctuations}, being a Central Limit Theorem (C.L.T.).
The study of the rate at which the probability of observing the discrete deviates from the expected limit decreases (roughly, exponentially fast) is the theme of the Large Deviations Principle.

Recently, the scientific community has given attention  to particle systems in random and non-homogeneous media, and several approaches have been developed in order to study the problem. In the papers \cite{fin, kk, nagy}, the authors considered random walks in a random environment, as for example the case where the environment is driven by an $\alpha$-subordinator. These works inspired a series of other papers in the context of particle systems, as \cite{fjl,fl,j,v}. The work in \cite{fl} was related to the hydrodynamic limit of exclusion processes driven by a general increasing function $W$, not necessarily a toss of an $\alpha$-subordinator. This work, in its hand, inspired the work \cite{fgn}, which dealt with the case $W$ being the distribution function of the Lebesgue measure plus a delta of Dirac measure, being the mass of the delta of Dirac dependent on the scale parameter. The model of \cite{fgn} can be described as follows. To each site of the discrete torus with $n$ sites, it is allowed to have at
most one particle. Each bond has a Poisson clock which is independent of the clocks on other sites. When the Poisson clock of a bond rings, the occupation at the vertices of this bond are interchanged. All the  Poisson clocks have parameter one, except one special clock, which has parameter given by $\alpha n^{-\beta}$, with $\alpha>0$ and $\beta\in[0,\infty]$. This ``slower" clock, makes the passage of particles across the corresponding bond more difficult, and for that reason that bond coined the name \emph{slow bond}.

In the scenario of \cite{fgn}, according to the value of $\beta$, three different limits for the time trajectory of the spatial density of particles were obtained. If $\beta\in [0,1)$ the limit is given by the weak solution of the periodic heat equation, meaning that the slow bond is not slow enough to originate any change in the continuum. If $\beta=1$, the limit is given by the weak solution of the heat equation with some Robin's boundary conditions representing the Fick's Law of passage of particles. And if $\beta\in{(1,\infty]}$, the limit is given by the weak solution of the heat equation with Neumann's boundary conditions, meaning that the slow bond in this regime of $\beta$ is slow enough to divide the space in the continuum.

Such dynamical phase transition (based on the strength of a single slow bond) is not limited to the hydrodynamic limit. In the ensuing papers \cite{fgn2,fgn3}, some other dynamical phase transitions were proved. In \cite{fgn2}, it was shown that the solutions of the three partial differential equations aforementioned are continuously related to a given boundary's parameter, indicating a dynamical phase transition also at the macroscopic level. In \cite{fgn3}, it was proved that the equilibrium fluctuations of the exclusion process with a slow bond evolving on an infinite volume, is also characterized by the same regimes of $\beta$. As before, in each case, namely for $\beta\in [0,1)$, $\beta=1$ or $\beta\in{(1,\infty]}$, the limit fluctuations of the system are driven by three Ornstein-Uhlenbeck processes. As a consequence of the density fluctuations, we have also obtained the corresponding phase transition for the current of particles through a fixed bond and for a
tagged particle.

In these notes  we make a synthesis of last results, all of them related to dynamic phase transitions that occur when the strength of a particular slow bond varies. We notice that the theme is not finished at all. There are a lot of particle systems to examine and different limits to prove. As an example, in the cited papers \cite{fgn,fgn2,fgn3}, the underlying particle systems are only of exclusion constrain and with symmetric dynamics. Therefore, one can exploit other dynamics and obtain other partial differential equations of physical interest. Moreover, even for the symmetric exclusion dynamics with a slow bond, the full scenario for the scaling limits is not closed yet: a Large Deviations Principle is still open. This is subject for future work.

 Here follows an outline of these notes. In Section 2 we present the  exclusion process with a slow bond. Section 3 is devoted to the scaling limits at the level of hydrodynamics. We present the hydrodynamic equations, the hydrodynamic limit and the phase transition for the corresponding partial differential equations.
In Section 4 we present the scaling limits at the level of fluctuations. We present the Ornstein-Uhlenbeck processes and the fluctuations of the density of particles. We finish in Section 5 with a description of the fluctuations of the current of particles and of a tagged particle.

\section{Exclusion processes}\label{s3}~

We are concerned with the study of dynamical phase transitions in particle systems with a single slow bond. Before discussing what we mean by a dynamical phase transition we describe our particle systems. We consider the simple exclusion process (SEP) with a single slow bond. Probabilistic speaking, the SEP  is a Markov process that we denote by $\{\eta_t:\, t\geq{0}\}$ and we consider it evolving on the state space $\Omega:=\{0,1\}^{\bb T_n}$, where $\bb T_n=\bb Z/n\bb Z$ is the one-dimensional discrete torus with $n$ points. A configuration of this Markov process is denoted by $\eta$ and it consists in a vector with $n$ components, each one taking the value $0$ or $1$. The physical interpretation is that whenever $\eta(x)=1$ we say that the site $x$ is occupied, otherwise it is empty.

The microscopic dynamics of this process can be informally described as follows. At each bond $\{x,x+1\}$ of $\bb{T}_n$, there is an exponential clock
of parameter $a^{n}_{x,x+1}$. When this clock rings, the value of $\eta$ at the vertices of this bond are interchanged. We choose the parameters of the clocks in all bonds equal to $1$, except at the bond $\{-1,0\}$, in such a way that the passage of particles across this bond is more difficult with respect to other bonds. For $\beta\in{[0,\infty]}$ and $\alpha>0$, we consider
   \begin{equation*}
a^{n}_{x,x+1}\;=\;\left\{\begin{array}{cl}
\alpha n^{-\beta}, &  \mbox{if}\,\,\,\,x=-1\,,\\
1, &\mbox{otherwise\,.}
\end{array}
\right.
\end{equation*}

This means that particles  cross all the bonds at rate $1$, except the bond $\{-1,0\}$, whose dynamics is slowed down as $\alpha n^{-\beta}$, with  $\alpha>0$ and $\beta\in{[0,\infty]}$, see the figure below.

\begin{figure}
\centering
\includegraphics{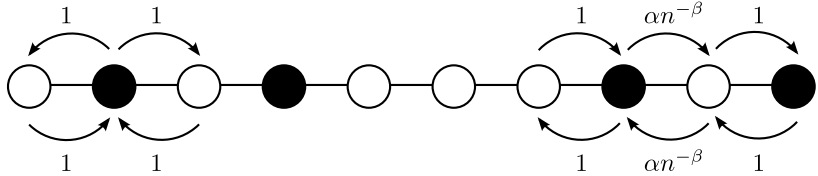}
% Fig1.eps: 0x0 pixel, 300dpi, 0.00x0.00 cm, bb=0 -1 337 67
\caption{SEP with a slow bond with vertices $\{-1,0\}$, whose jump rates are given by $\alpha n^{-\beta}$. Black balls represent occupied sites.}
\label{fig:1}
\end{figure}

The dynamics described above can be characterized via the infinitesimal generator, which we denote by $\mathcal{L}_{n}$ and is given on functions $f:\Omega\rightarrow \bb{R}$ as
\begin{equation*}
\mathcal{L}_{n}f(\eta)=\sum_{x\in \bb T_n}\,a^{n}_{x,x+1}\,\big[f(\eta^{x,x+1})-f(\eta)\big]\,,
\end{equation*}
where $\eta^{x,x+1}$ is the configuration obtained from $\eta$ by exchanging the occupation variables $\eta(x)$ and $\eta(x+1)$, namely
\begin{equation*}
\eta^{x,x+1}(y)=\left\{\begin{array}{cl}
\eta(x+1),& \mbox{if}\,\,\, y=x\,,\\
\eta(x),& \mbox{if} \,\,\,y=x+1\,,\\
\eta(y),& \mbox{otherwise}\,.
\end{array}
\right.
\end{equation*}

Let $\rho\in{[0,1]}$ and denote the Bernoulli product measure, defined in  $\Omega$ and with parameter $\rho$, by $$\nu^n_\rho\{\eta\in\Omega: \,\,\eta(x)=1,\;\;\mbox{for any}\;\; x\in A\}=\rho^{\# A},$$ for all set $A\subset \bb T_n$. Here $\# A$ denotes the cardinality of the set $A$. It is well known that the  measures  $\nu^n_\rho$  are invariant for the dynamics   introduced above.  Moreover,  these measures are also
reversible.

 The trajectories of the Markov process $\{\eta_t : t\ge 0\}$ live on the space $\mc D(\bb R_+, \Omega)$, that is, the path space of
c\`adl\`ag trajectories with values in $\Omega$. For a
measure $\mu_n$ on $\Omega$, we denote by  $\bb P_{\mu_n}$ the
probability measure on $\mc D(\bb R_+, \Omega)$ induced by $\mu_n$ and $\{\eta_t : t\ge 0\}$; and we denote by  $\bb E_{\mu_n}$
expectation with respect to $\bb P_{\mu_n}$.

We notice that we do not index the Markov process, the generator nor the measures, in $\beta$ or $\alpha$ for simplicity of notation.

\section{Hydrodynamical phase transition}

The study of the hydrodynamical behavior consists in the analysis of the time evolution of the density of particles. For that purpose  we introduce the empirical measure process as follows.

 For $t\in{[0,T]}$, let $\pi^{n}_t(\eta,du):=\pi^n(\eta_t,du)\in \mc M$ be defined as
\begin{equation*}%\label{f01}
\pi^{n}(\eta_t,du) \;=\; \pfrac{1}{n} \sum _{x\in \bb T_n} \eta_t (x)\,
\delta_{x/n}(du)\,,
\end{equation*}
where $\delta_y$ is the Dirac measure concentrated on $y\in \bb T$. Above,  $\bb T$ denotes the one-dimensional torus and  $\mc M$ denotes the space of positive measures on $\bb T$ with total
mass bounded by one, endowed with the weak topology.

The hydrodynamic limit can be stated as follows. If we assume a L.L.N. for $\{\pi^n_0\}_{n\in\bb N}$ to a limit $\rho_0(u)du$ under the initial distribution of the system, then at any time $t>0$ the L.L.N. holds for $\{\pi^n_t\}_{n\in\bb N}$ to a limit $\rho(t,u)du$ under the corresponding distribution of the system at time $t$. Moreover, the density $\rho(t,u)$ evolves according to a partial differential equation - the hydrodynamic equation.
For this model, depending on the  range of the parameter $\beta$, we obtain different hydrodynamic equations for the underlying particle system.

In the next section we describe the hydrodynamic equations we obtained and we precise in which sense $\rho(t,u)$ is a solution to those equations.

\subsection{Hydrodynamic equations}~

We start by describing the hydrodynamic equations that govern the evolution of the density of particles for the models introduced above.
Depending on the range of the parameter $\beta$ we obtain hydrodynamic equations which have different behavior. More precisely, we always obtain the heat equation but with different boundary conditions. The first hydrodynamic equation is the heat equation with periodic boundary conditions, namely:
 \begin{equation}   \label{he}
\left\{
\begin{array}{ll}
 \partial_t \rho(t,u) \; =\; \Delta \rho(t,u)\,, \qquad& t \geq 0,\, u\in \bb T\,,\\
  \rho(0,u) \;=\; \rho_0(u), &u \in \bb T\,.
\end{array}
\right.
\end{equation}

In the hydrodynamic limit scenario, we obtain $\rho(t,u)$ as  a weak solution of the corresponding hydrodynamic equation. To  make this notion precise, we introduce the following definition:

 \begin{definition}\label{def edp 1}
 Let $\rho_0:\bb T\to [0,1]$ be a measurable function.
   We say that $\rho:[0,T]\times \bb T\to[0,1]$ is a weak solution of the heat equation with periodic boundary conditions given in \eqref{he} if $\rho$ is measurable and, for any  $t\in{[0,T]}$ and any $H\in C^{1,2}([0,T]\times\mathbb{T}) $,
\begin{equation}\label{eqint1}
\begin{split}
\int_{\bb T}\rho(t,u)H(t,u) du-&\int_{\bb T}\rho(0,u)H(0,u) du \\&-\int_0^t\int_{\bb T}\rho(s,u)\Big(\partial_s H(s,u)+\Delta  H(s,u)\Big) du\;ds\;=\;0\,.
\end{split}
\end{equation}
\end{definition}
Above and in the sequel the space $C^{1,2}([0,T]\times\mathbb{T})$ is the space of  real valued functions defined on $[0,T]\times\mathbb{T}$
of class $C^1$ in time and $C^2$ in space.

The second equation we consider is the heat equation with a type of Robin's boundary conditions, that is:
\begin{equation}\label{her}
\left\{
\begin{array}{ll}
 \partial_t \rho(t,u) \; =\; \Delta \rho(t,u)\,, &t \geq 0,\, u\in (0,1)\,,\\
 \partial_u \rho(t,0) \; =\;\partial_u \rho(t,1)= \alpha(\rho(t,0)-\rho(t,1))\,, \qquad &t \geq 0,\\
 \rho(0,u) \;=\; \rho_0(u), &u \in (0,1)\,.
\end{array}
\right.
\end{equation}

To introduce the notion of weak solution of this equation we need to recall  the notion of Sobolev's spaces.

\begin{definition}\label{Sobolevdefinition}
Let  $\mc H^1$ be the set of all locally summable functions $\zeta: (0,1)\to\bb R$ such that
there exists a function $\p_u\zeta\in L^2(0,1)$ satisfying
\begin{equation*}
 \int_{\bb T}\partial_uG(u)\zeta(u)\;du=\,-\int_{\bb T}G(u)\partial_u\zeta(u)\;du\,,
\end{equation*}
for all $G\in C^{\infty}(0,1)$ with compact support.
Let $L^2(0,T;\mc H^1)$ be the space of
 all measurable functions
$\xi:[0,T]\to \mc H^1$ such that
\begin{equation*}
\Vert\xi \Vert_{L^2(0,T;\mc H^1)}^2 \,
:=\,\int_0^T \Big(\Vert \xi\Vert_{L^2[0,1]}^2+\Vert\partial_u\xi\Vert_{L^2[0,1]}^2\Big)\,dt\,<\,\infty\,.
\end{equation*}
\end{definition}
Above $\Vert \cdot \Vert_{L^2[0,1]}$ denotes the $L^2$-norm in $[0,1]$.

\begin{definition}\label{heat equation Robin}
Let $\rho_0:\bb T\to [0,1]$ be a measurable function. We say that $\rho:[0,T]\times \bb T\to[0,1]$  is a weak solution of the heat equation with Robin's boundary conditions given in \eqref{her}
if $\rho\in L^2(0,T;\mathcal{H}^1)$ and for  all $t\in [0,T]$ and for all  $H\in \A$,
\begin{equation}\label{eqint2}
\begin{split}
\int_{\bb T}\rho(t,u)H(t,u) du&-\int_{\bb T}\rho(0,u)H(0,u) du \\&-\int_0^t\int_{\bb T}\rho(s,u)\Big(\partial_s H(s,u)+\Delta  H(s,u)\Big) du\;ds\\&-\!\!\int_0^t\!\!\!(\rho_s(0)\partial_uH_s(0)-\rho_s(1)\partial_uH_s(1))\,ds\\
&+\int_0^t \alpha(\rho_s(0)-\rho_s(1))(H_s(0)-H_s(1))\,ds=0\,.
\end{split}
\end{equation}
\end{definition}
The last equation we consider is the heat equation with Neumann's boundary conditions given by:
 \begin{equation}\label{hen}
\left\{
\begin{array}{ll}
 \partial_t \rho(t,u) \; =\; \Delta \rho(t,u)\,, &t \geq 0,\, u\in (0,1)\,,\\
 \partial_u \rho(t,0) \; =\;\partial_u \rho(t,1)= 0\,, \qquad &t \geq 0\,,\\
 \rho(0,u) \;=\; \rho_0(u), &u \in (0,1)\,.
\end{array}
\right.
\end{equation}

\begin{definition}\label{heat equation Neumann}
Let $\rho_0:\bb T\to [0,1]$ be a measurable function. We say that $\rho:[0,T]\times \bb T\to[0,1]$  is a weak solution of the heat equation with Neumann's boundary conditions if $\rho\in L^2(0,T;\mathcal{H}^1)$ and for all $t\in [0,T]$ and for all  $H\in \A$,
\begin{equation}\label{eqint3}
\begin{split}
\int_{\bb T}\rho(t,u)H(t,u) du&-\int_{\bb T}\rho(0,u)H(0,u) du \\&-\int_0^t\int_{\bb T}\rho(s,u)\Big(\partial_s H(s,u)+\Delta  H(s,u)\Big) du\;ds\\&-\!\!\int_0^t\!\!\!(\rho_s(0)\partial_uH_s(0)-\rho_s(1)\partial_uH_s(1))\,ds=0\,.
\end{split}
\end{equation}
\end{definition}

Our argument to prove the hydrodynamic limit is standard in the theory of stochastic processes and goes through a tightness argument for $\{\pi^n_t\}_{n\in\bb N}$, which means relatively compactness of $\{\pi^n_t\}_{n\in\bb N}$. Therefore, there exists a limit point. To have uniqueness of the limit point of $\{\pi^n_t\}_{n\in\bb N}$ it is sufficient to prove uniqueness of the weak solution of the corresponding  hydrodynamic equation. Then, it follows the convergence of the whole sequence $\{\pi^n_t\}_{n\in\bb N}$ to the unique limit point.  For tightness issues we refer the reader to \cite{fgn} and the uniqueness of the weak solution is stated below.

\begin{proposition}
Let $\rho_0:\bb T\to [0,1]$ be a measurable function. There exists a unique weak solution of the heat equation with periodic boundary conditions given in \eqref{he} and  a unique weak solution of the heat equation with Neumann's boundary conditions given in \eqref{hen}. Moreover, for each $\alpha>0$, there exists a unique weak solution of the heat equation with Robin's boundary conditions given in \eqref{her}.
\end{proposition}

\subsection{Hydrodynamic limit}~

Returning to our discussion on the validity of the hydrodynamic limit, we introduce the set of initial measures for which we deduce the result.

\begin{definition} \label{def associated measures}
Let $\rho_0:\bb T\to [0,1]$ be a measurable function.
A sequence of probability measures $\{\mu_n\}_{n\in\bb N}$ on $\Omega$ is
said to be associated to a profile $\rho_0 :\bb T \to [0,1]$ if, for every $\delta>0$ and every continuous function  $H:\bb T\to \bb R$, it holds that
\begin{equation}\label{associated}
\lim_{n\to\infty}
\mu_n \Big\{ \eta:\, \Big\vert \pfrac 1n \sum_{x\in\bb T_n} H(\pfrac{x}{n})\, \eta(x)
- \int_{\bb{T}} H(u)\, \rho_0(u) du \Big\vert > \delta \Big\}\;=\; 0\,.
\end{equation}
\end{definition}

One could ask about the existence of a measure associated to the profile $\rho_0:\bb T\to [0,1]$. For instance, we can consider a Bernoulli product measure in $\Omega$ with marginal at $\eta(x)$ given by $\mu_n\{\eta\in \Omega: \eta(x)=1\}=\rho_0(x/n)$.

For these processes we obtained in \cite{fgn,fgn2} that:

\begin{theorem} \label{th:hlrm} [L.L.N. for the density of particles]
Fix $\beta\in [0,\infty]$ and $\rho_0: \mathbb{T} \to [0,1]$ a measurable function. Let $\{\mu_n\}_{n\in\bb N}$ be
a sequence of probability measures  on $\Omega$ associated to $\rho_0$. Then, for any $t\in [0,T]$, for every $\delta>0$ and every continuous function  $H:\bb T\to \bb R$:
\begin{equation*}
\lim_{n\to\infty}
\mathbb{P}_{\mu_n} \Big\{\eta_. : \, \Big\vert \pfrac{1}{n} \sum_{x\in\mathbb{T}_n}
H\big(\pfrac{x}{n}\big)\, \eta_t(x) - \int_{\bb T}H(u)\rho(t,u)du \Big\vert
> \delta \Big\}\;=\; 0\,,
\end{equation*}
 where:
\begin{itemize}
\item
for $\beta\in[0,1)$, $\rho(t,\cdot)$ is the unique weak solution of \eqref{he};

\item
for $\beta=1$, $\rho(t,\cdot)$ is the unique weak solution of \eqref{her};

\item
 for $\beta\in(1,\infty]$, $\rho(t,\cdot)$ is the unique weak solution of  \eqref{hen}.
\end{itemize}
All equations have the same initial condition $\rho_0:\bb T\to [0,1]$.
\end{theorem}

\subsection{Phase transition for the hydrodynamic equations}\label{pt}~

A puzzling question is whether there is a similar phase transition as described above, but at the macroscopic level. More precisely, does the unique weak solution of the heat equation with Robin's boundary conditions, that we denote by $\rho^\alpha$, converge in any sense to the weak solution of the heat equation with periodic boundary conditions or to the weak solution of the heat equation with Neumann's boundary conditions?  In \cite{fgn2} we gave an affirmative answer to this question. We proved that $\rho^\alpha$ converges to the unique weak solution of the heat equation with Neumann's
boundary conditions, when $\alpha$ goes to zero and to the unique weak solution of the heat equation with periodic
boundary conditions, when $\alpha$ goes to infinity. This is the content of the next theorem.

This result is concerned only with the partial differential equations, having at principle nothing to do with the underlying particle systems. Nevertheless, our approach of proof is based on energy estimates coming from these particle systems.

\begin{theorem}\label{pdePT} [Phase transition for the heat equation with Robin's boundary conditions]
 For $\alpha>0$, let $\rho^\alpha:[0,T]\times[0,1]\to [0,1]$ be the unique weak solution of the heat equation with Robin's boundary conditions:
 \begin{equation*}
\left\{
\begin{array}{ll}
 \partial_t \rho^\alpha(t,u) \; =\; \Delta \rho^\alpha(t,u)\,, &t \geq 0,\, u\in (0,1)\,,\\
 \partial_u \rho^\alpha(t,0) \; =\;\partial_u \rho^\alpha(t,1)= \alpha(\rho^\alpha(t,0)-\rho^\alpha(t,1))\,,  &t \geq 0\,,\\
 \rho^\alpha(0,u) \;=\; \rho_0(u), &u \in (0,1)\,.
\end{array}
\right.
\end{equation*}
Then, $\lim_{\alpha\to 0} \rho^\alpha \; = \;  \rho^0,$ in $L^2([0,T]\times [0,1])$, where
$\rho^0:[0,T]\times[0,1]\to [0,1]$ is the unique weak solution of the heat equation with  Neumann's boundary  conditions
 \begin{equation*}
\begin{cases}
 \partial_t \rho^0(t,u) \; =\; \Delta \rho^0(t,u)\,,&t \geq 0,\, u\in (0,1)\,,\\
 \partial_u \rho^0(t,0) \; =\;\partial_u \rho^0(t,1)= 0\,,&t \geq 0\,,\\
  \rho^0(0,u) \;=\; \rho_0(u)\,, &u \in (0,1)\,\\
\end{cases}
\end{equation*}
and $\lim_{\alpha\to \infty} \rho^\alpha \; = \;  \rho^\infty,$
in $L^2([0,T]\times [0,1])$,
 where $\rho^\infty:[0,T]\times[0,1]\to [0,1]$ is the unique weak solution of the heat equation with
periodic boundary conditions
\begin{equation*}
\begin{cases}
 \partial_t \rho^\infty(t,u) \; =\; \Delta \rho^\infty(t,u)\,,\qquad & t \geq 0,\, u\in \mathbb T\,,\\
\rho^\infty(0,u) \;=\; \rho_0(u)\,, &u \in \mathbb T\,.
\end{cases}
\end{equation*}

\end{theorem}

\section{Equilibrium fluctuations}~

Above we obtained a L.L.N. for the empirical measure considering the process starting from a measure which is associated to a profile $\rho_0:\bb T\to [0,1]$. The natural question that follows is: what are the fluctuations around this ``mean" profile? Do we have a C.L.T. for the density of particles? Under what set of initial measures? In the next lines we answer this question for a particular set of initial distributions, namely for the invariant measures $\nu^n_\rho$. In case of non-invariant measures the problem is still open.

In this case we consider the process evolving on $\bb Z$, being its state space $\{0,1\}^{\bb Z}$. To define properly our results, we fix $\rho\in[0,1]$, and  we introduce the density fluctuation field as follows.
 For $t\in{[0,T]}$, let
\begin{equation*}
\mc Y^{n}_t(\eta,du) \;=\; \sqrt n \pi^n_{tn^2}(\eta,du)-E_{\nu^n_\rho}[\sqrt n \pi^n_{tn^2}(\eta,du)],
\end{equation*}
where $x$ runs through $\bb Z$ in the definition of $\pi^n_t(\eta,du)$ and $E_{\nu^n_\rho}$ denotes expectation with respect to $\nu^n_\rho$.
Then, for any function $H:\bb R\to\bb R$ we have that
\begin{equation*}
\int_{\bb R}H(u)\mc Y^{n}_t(\eta,du) \;=\; \frac{ 1}{\sqrt n} \sum_{x\in{\bb Z}}H\Big(\frac{x}{n}\Big)[\eta_{tn^2}(x)-\rho].
\end{equation*}

By computing the characteristic function of $\mc{Y}^n_0$, we obtain that  $\{\mathcal{Y}_{0}^n\}_{ n\in{\bb N}}$ converges as $n$ goes to $\infty$ to a mean zero gaussian process $\mc Y_0$. More precisely, for any $H$, $\mc Y_0(H)$ is a gaussian random variable with mean zero and variance given by $$\rho(1-\rho)\int_\bb R (H(x))^2dx.$$ Next, we are going to  characterize the stochastic partial differential equations governing the evolution of the limit points of $\{\mc Y_t^n\}_{n\in{\mathbb{N}}}$.

\subsection{Ornstein-Uhlenbeck processes}~

In order to properly write down the stochastic partial differential equations that we deal with, we need to introduce different sets of test functions and
two type of operators defined on these spaces.

\begin{definition}
Define $\mc S(\bb R\backslash \{0\})$ as the space of functions $H\in C^\infty(\bb R\backslash\{0\})$, that are continuous from the right at $x=0$, for which
\begin{equation*}
 \Vert H \Vert_{k,\ell}\;:=\;\sup_{x\in \bb R\backslash{\{0\}}}|(1+|x|^\ell)
\,H^{(k)}(x)|\;<\;\infty\,,
\end{equation*}
\noindent for all integers $k,\ell\geq 0$, and $H^{(k)}(0^-)=H^{(k)}(0^+)$, for all $k$ integer, $k\geq 1$.
\begin{itemize}
\item For $\beta\in[0,1)$, let $\mc S_\beta(\bb R)$ be the subset of $\mc S(\bb R\backslash \{0\})$
composed of functions $H$ satisfying $H(0^-)=H(0^+)\,.$

\item For $\beta=1$,  let $\mc S_\beta(\bb R)$ as the subset of $\mc S(\bb R\backslash \{0\})$
composed of functions $H$ satisfying $H^{(1)}(0^+)\;=\; H^{(1)}(0^-)\;=\; \alpha ( H(0^+)-H(0^-))\,.$

\item For $\beta\in(1,+\infty]$,  let $\mc S_\beta(\bb R)$ be the subset of $\mc S(\bb R\backslash \{0\})$
composed of functions $H$ satisfying $H^{(1)}(0^+)\;=\; H^{(1)}(0^-)\;=\; 0\,.$
\end{itemize}
\end{definition}

Above and in the sequel, $H^{(k)}(\cdot)$ represents the $k$-th derivative of the function $H$ and $H(0^{+})$ (resp. $H(0^-)$) denotes the limit of $H$ from the right (resp. left) of $0$.

\begin{definition}
For $\beta\in{[0,\infty]},$ we define the operators  $\Delta_\beta,\nabla_\beta : \mc S_\beta(\bb R)\rightarrow \mc S(\bb R)$ by
\begin{equation*}
\nabla_\beta H(u)\;=\;\left\{\begin{array}{cl}
H^{(1)}(u), &  \mbox{if}\,\,\,\,u\neq 0\,,\\
H^{(1)}(0^+), &\mbox{if}\,\,\,\,u=0\,,
\end{array}
\right.
\end{equation*}
and
\begin{equation*}
\Delta_\beta H(u)\;=\;\left\{\begin{array}{cl}
H^{(2)}(u), &  \mbox{if}\,\,\,\,u\neq 0\,,\\
H^{(2)}(0^+), &\mbox{if}\,\,\,\,u=0\,,
\end{array}
\right.
\end{equation*}
\end{definition}
which are essentially the usual derivative and the usual second derivative, but defined in the domains $\mc S_{\beta}(\bb R)$. We have the following uniqueness result  which is a key point in our approach.

Denote by $T^\beta_{t}$ the semigroup corresponding to the partial differential equations \eqref{he}, \eqref{her} or  \eqref{hen}, if $\beta\in [0,1)$, if $\beta=1$ or if $\beta\in(1,\infty]$, respectively.
\begin{proposition}\label{pp1}
For each $\beta\in[0,\infty]$ and $\alpha>0$, there exists an unique random element $\mc Y_\cdot$ taking values in the space $C([0,T],\mathcal{S}'_{\beta}(\bb R))$ such
that:
\begin{itemize}
\item[i)] For every function $H \in \mathcal{S}_{\beta}(\bb R)$, $\mc M_t(H)$ and $\mc N_t(H)$ given by
\begin{equation}\label{lf1}
\begin{split}
&\mc M_t(H)= \mc Y_t(H) -\mc Y_0(H) -  \int_0^t \mc Y_s(\Delta_\beta H)ds\,,\\
&\mc N_t(H)=\big(\mc M_t(H)\big)^2 - 2\chi(\rho) \; t\,\|\nabla_\beta H\|_{2,\beta}^2
\end{split}
\end{equation}
are $\mc F_t$-martingales, where  $\mc F_t:=\sigma(\mc Y_s(H); s\leq t,  H \in \mathcal{S}_{\beta}(\bb R))$, for $t\in{[0,T]}$.
\item[ii)] $\mc Y_0$ is a mean zero gaussian field with covariance given on $G,H\in{\mathcal{S}_{\beta}(\mathbb{R})}$ as
\begin{equation}\label{eq:covar1}
\mathbb{E}\big[ \mc Y_0(G) \mc Y_0(H)\big] =  \chi(\rho)\int_{\mathbb{R}} G(u) H(u) du\,.
\end{equation}
\end{itemize}
Moreover, for each $H\in\mc S_\beta(\bb R)$, the stochastic process $\{\Y_t(H)\,;\,t\geq 0\}$ is  gaussian, being the
distribution of $\Y_t(H)$  conditionally to
$\mc F_s$, for $s<t$, gaussian of mean $\Y_s(T^\beta_{t-s}H)$ and variance $\int_0^{t-s}\Vert \nabla_\beta T^\beta_{r}
H\Vert^2_{2,\beta}\,dr$.
\end{proposition}
Above and in the sequel $\mathcal{S}'_{\beta}(\mathbb{R})$ denotes the space of bounded linear functionals $f:
\mathcal{S}_{\beta}(\mathbb{R})\rightarrow{\bb R}$ and  $\mathcal{D}([0,T],\mathcal{S}'_{\beta}(\mathbb{R}))$
(resp. $C([0,T],\mathcal{S}'_{\beta}(\mathbb{R}))$)
 is the space of  c\`adl\`ag (resp. continuous) $\mathcal{S}'_{\beta}(\mathbb{R})$ valued functions endowed with the Skohorod topology.
Also $\|H\|_{2,\beta}^2=\|H\|_2^2+(H(0))^2\textbf{1}_{\{\beta=1\}},$ where $\|\cdot\|_{2}$ denotes the $L^2$-norm in $\bb R$.
 We call to $\mc Y_\cdot$ the generalized Ornstein-Uhlenbeck process of
characteristic operators $\Delta_\beta$ and $\nabla_\beta $ and it is
the formal solution of the following equation $$d\mathcal{Y}_t=\Delta_\beta \mathcal{Y}_tdt+\sqrt{2\chi(\rho)}\nabla_\beta d\mc{W}_t\,,$$ where $\mc{W}_t$ is a space-time white noise of unit variance.

\subsection{Central Limit Theorem}~

We are in position to state the equilibrium fluctuations for the density of particles. Notice that our initial distribution is $\nu^n_{\rho}$, an invariant measure.
\begin{theorem}[C.L.T. for the density of particles]\label{flu1} The sequence of processes $\{\mathcal{Y}_{t}^n\}_{ n\in{\bb N}}$ converges in distribution, as $n$ goes to $\infty$, with respect to the
Skorohod topology
of $\mathcal{D}([0,T],\mathcal{S}'_{\beta}(\bb R))$ to a gaussian process $\mathcal{Y}_t$
in $C([0,T],\mathcal{S}'_{\beta}(\bb R))$, which is the formal solution of the Ornstein-Uhlenbeck equation given by
\begin{equation}\label{eq Ou}
d\mathcal{Y}_t=\Delta_\beta \mathcal{Y}_tdt+\sqrt{2\chi(\rho)}\nabla_\beta d\mc{W}_t\,.
\end{equation}
\end{theorem}

\section{Current and Tagged particle fluctuations}~
In this section we are still restricted to the invariant state $\nu_{\rho}^n$ and for that purpose we fix a density $\rho$ from now on up to the rest of these notes.
\subsection{The current}~
Now, we introduce the notion of current of particles through a fixed bond $\{x,x+1\}$. For a bond $e_x:=\{x,x+1\}$, denote
by ${J}^n_{e_x}(t)$ the current of particles over the bond $e_x$, that is ${J}^n_{e_x}(t)$ counts the total number of jumps from the
site $x$ to the site $x+1$ minus
the total number of jumps from the site $x+1$ to the site $x$ in the time interval $[0,tn^2]$, see the figure below.
\begin{figure}
\centering
\includegraphics{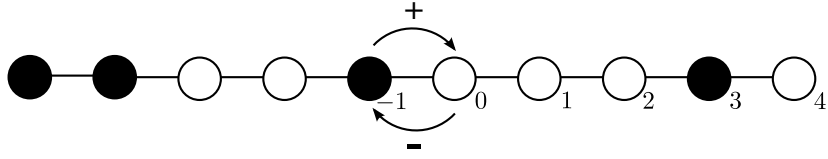}
% Fig1.eps: 0x0 pixel, 300dpi, 0.00x0.00 cm, bb=0 -1 337 67
\caption{Current at the bond $\{-1,0\}$ of the SEP with a slow bond. Every time a particle jumps from $-1$ to $0$ ($0$ to $-1$) the current increases (decreases) by one.}
\label{fig:2}
\end{figure}
More generally, to each macroscopic point $u\in{\mathbb{R}}$ we can define the current through its associated microscopic bond of vertices
$\{\lfloor un \rfloor -1, \lfloor un \rfloor\}$, as ${J}^n_{u}(t):= {J}^n_{e_{\lfloor un \rfloor -1}}(t)\,.$ Here $\lfloor un\rfloor$ denotes  the biggest integer smaller or equal to $un$. As a consequence of the C.L.T. for the density of particles, namely of Theorem \ref{flu1}, it is simple to derive the C.L.T. for the current of particles which we enounce as follows.

\begin{theorem}[C.L.T. for the current of particles] \label{current clt}
Under $\bb P_{\nu_\rho^n}$, for every $t\geq{0}$ and every $u\in\bb R$, $$\frac{{J}^n_u(t)}{\sqrt n}\xrightarrow[{n\rightarrow{\infty}}]\,{J}_u({t})$$
in the sense of finite-dimensional distributions, where ${J}_u({t})$  is  a gaussian process with  mean zero and variance given by
\medskip

$\bullet$ for $\beta\in{[0,1)}$, $\mathbb{E}_{\nu_\rho^n}[({J}_u(t))^2]=2\chi(\rho)\sqrt{\frac{t}{\pi}}$,
that is $J_u(t)$ is a fractional Brownian Motion of Hurst exponent $1/4$;

\quad

$\bullet$ for $\beta=1$, $\mathbb{E}_{\nu_\rho^n}[({J}_u(t))^2]=2\chi(\rho)\Big(\sqrt{\frac{t}{\pi}}+\frac{\Phi_{2t}(2u+4\alpha
t)\,e^{4\alpha u+4\alpha^2t}-\Phi_{2t}(2u)}{2\alpha}\Big)$;

\quad

$\bullet$ for  $\beta\in (1,{+\infty}]$, $\mathbb{E}_{\nu_\rho^n}[({J}_u(t))^2]=2\chi(\rho)\Big(\sqrt{\frac{t}{\pi}}\Big[1-e^{-u^2/t}\Big]+2u\,\Phi_{2t}(2u)\Big),$

where \begin{equation*}
\Phi_{2t}(x):=\int^{{+\infty}}_x\frac{e^{-u^2/{4t}}}{\sqrt{4\pi t}}du \,.
\end{equation*}
\end{theorem}
It worth to remark the variance at $u=0$, corresponding to the current of particles through the slow bond $e_{-1}$. If $\beta\in [0,1)$,
the variance corresponds to the one of a fractional Brownian Motion of Hurst exponent $1/4$. If $\beta\in (1, {\infty}]$, the
variance equals to zero as expected. This is a consequence of having Neumann's boundary conditions at $x=0$ which turns it into an isolated boundary. And for $\beta=1$, we
obtain a family of gaussian processes indexed in $\alpha$ interpolating the two aforementioned processes.

\begin{corollary} \label{limit robin current}

For $\beta=1$, denote the limit, as $n\to\infty$, of $J^n_u(t)/\sqrt n$ by
${J}^\alpha_u({t}).$

 Then for every $t\geq{0}$ and every $u\in\bb R$, $${J}^\alpha_u({t})\xrightarrow[{\alpha\rightarrow{+\infty}}]\,{J}^\infty_u({t}),$$
 where $J^\infty_u(t)$ is the fractional  Brownian Motion
with Hurst exponent $1/4$
and $${J}^\alpha_u({t})\xrightarrow[{\alpha\rightarrow{0}}]\,{J}^0_u({t})\,,$$
where $J^0_u(t)$ is the mean zero gaussian process with variance given by \\ $\mathbb{E}_{\nu_\rho^n}[({J}_u(t))^2]=2\chi(\rho)\Big(\sqrt{\frac{t}{\pi}}\Big[1-e^{-u^2/t}\Big]+2u\,\Phi_{2t}(2u)\Big).$

 The convergence is in the sense of finite dimensional distributions.
 \end{corollary}

\subsection{Tagged particle fluctuations}~

Our last goal is to present the asymptotic behavior of a tagged particle in the system. The dynamic of this tagged particle is no longer Markovian, since its behavior is influenced by the presence of other particles in the system. Nevertheless, we can relate the position of the tagged particle with the current and the density of particles, and from the previous results we obtain information about the behavior of this particle.

Suppose to start the system from a configuration with a particle at the site $\lfloor un\rfloor$ and in all other sites suppose that the configuration is distributed according to $\nu^n_\rho$. In other words, this means that we consider the Markov process $\{\eta_t:t\geq 0\}$ starting from the measure $\nu^n_{\rho}$ conditioned to have a particle at the site $\lfloor un\rfloor$, that we denote by $\nu_\rho^{n,u}$. That is, $\nu_{\rho}^{u,n}(\cdot):=\nu^n_\rho(\,\cdot\,|\eta_{tn^2}(\lfloor un\rfloor) =1)$.

\begin{figure}
\centering
\includegraphics{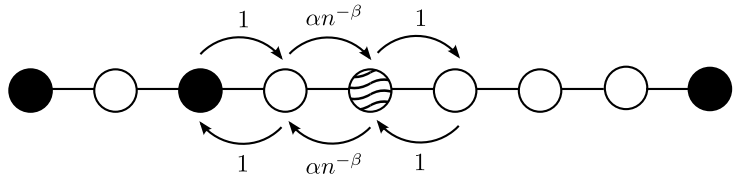}
% Fig1.eps: 0x0 pixel, 300dpi, 0.00x0.00 cm, bb=0 -1 337 67
\caption{The tagged particle of the SEP with a slow bond. At initial time, the tagged particle is at the site $0$.}
\label{fig:3}
\end{figure}

We notice that the previous results were obtained considering the process starting from $\nu^n_\rho$. In order to be able to use them, we couple the process starting from $\nu_\rho^{n,u}$ and starting from $\nu^n_\rho$, in such a way that both
processes  differ at most by one site at any given time. This allow us to derive the same statements of Theorems \ref{flu1} and \ref{current clt} for the starting measure  $\nu_\rho^{n,u}$.

  Now, let $X^n_u(t)$ be the position at the time $tn^2$ of a tagged particle initial at the site $\lfloor un\rfloor$.
Since our study is restricted to the one dimensional setting, particles do preserve their order, and it is simple to check that
\begin{equation*}\label{relation tp cur dp}
\{X^n_u(t)\geq{k}\}=\Big\{ {J}^n_{u}(t)\geq{\sum_{x=\lfloor un \rfloor}^{\lfloor un \rfloor+k-1}\eta_{tn^2}(x)}\Big\}.
\end{equation*}

We explain briefly how to get the previous equality. Suppose for simplicity that $u=0$, so that we start the system with the tagged particle at the origin. If this particle is, at time $tn^2$, at the right hand side of $n$, then all the particles that jumped from  $-1$ to $0$ and did not jump backwards, are somewhere at the sites $\{0,1,\ldots,X^n_u(t)\}$. It follows that the current through the bond $\{-1,0\}$ has to be greater or equal than the density of particles in $\{0,\ldots,n\}$. Reasoning similarly, we get the equality between those events.

 Finally, last relation together with Theorem \ref{current clt}, implies the following result.
 \begin{theorem}[C.L.T. for a tagged particle] \label{tagged clt}
Under $\bb{P}_{\nu^u_\rho}$, for all $\beta\in{[0,\infty]}$,  every $u\in\bb R$ and $t\geq{0}$
\begin{equation*}
\frac{{X}^n_u(t)}{\sqrt{n}}\xrightarrow[{n\rightarrow{+\infty}}]\,X_u(t)
\end{equation*}
in the sense of finite-dimensional distributions, where $X_u(t)=J_u(t)/\rho$ in law and $J_u(t)$ is the same as in Theorem
\ref{current clt}. In particular, the variance of the process $X_u(t)$ is given by
\medskip

$\bullet$ for $\beta\in{[0,1)}$, $\mathbb{E}_{\nu_\rho^n}[({X}_u(t))^2]=2\cfrac{\chi(\rho)}{\rho^2}\sqrt{\frac{t}{\pi}}$,
that is $X_u(t)$ is a fractional Brownian Motion of Hurst exponent $1/4$;

\quad

$\bullet$ for $\beta=1$, $\mathbb{E}_{\nu_\rho^n}[({X}_u(t))^2]=2\cfrac{\chi(\rho)}{\rho^2}\Big(\sqrt{\frac{t}{\pi}}+\frac{\Phi_{2t}(2u+4\alpha
t)\,e^{4\alpha u+4\alpha^2t}}{2\alpha}\Big)$;

\quad

$\bullet$ for  $\beta\in (1,{+\infty}]$, $\mathbb{E}_{\nu_\rho^n}[({X}_u(t))^2]=2\cfrac{\chi(\rho)}{\rho^2}\Big(\sqrt{\frac{t}{\pi}}\Big[1-e^{-u^2/t}\Big]+2u\,\Phi_{2t}(2u)\Big).$

\end{theorem}

\section*{Acknowledgements}~~~
The authors thank the great hospitality of CMAT (Portugal), IMPA and PUC (Rio de Janeiro).

A.N. thanks Cnpq (Brazil) for support through
the research project ``Mecânica estatística fora do equilíbrio para sistemas estocásticos" Universal n. 479514/2011-9.

P.G. thanks FCT (Portugal) for support through the research
project ``Non-Equilibrium Statistical Physics" PTDC/MAT/109844/2009.
P.G. thanks the Research Centre of Mathematics of the University of
Minho, for the financial support provided by ``FEDER" through the
``Programa Operacional Factores de Competitividade  COMPETE" and by
FCT through the research project PEst-C/MAT/UI0013/2011.

T.F. was supported through a grant ``BOLSISTA DA CAPES - Brasília/Bra-sil" provided by CAPES (Brazil).

\bibliographystyle{amsplain}

\end{document}